\documentstyle[11pt,amstex,amssymb]{amsart}\textheight=19.8truecm\textwidth=13.7truecm
\begin{document} 

\vskip 30pt
\centerline{\Large\bf Symmetric norms and reverse inequalities to}
\vskip 20pt 
\centerline{\Large\bf   Davis and Hansen-Pedersen characterizations}
\vskip 20pt
\centerline{\Large\bf of operator convexity}
\vskip 30pt
\centerline{Jean-Christophe Bourin}
\vskip 10pt
 \centerline{E-mail: bourinjc@@club-internet.fr}
 \vskip 10pt
 \centerline{\it Dedicated to Fran\c coise Piquard, respectfully and affectionately}
\vskip 25pt
\noindent
{\small {\bf Abstract.} Let $A$, $B$, $Z$ be $n$-by-$n$ matrices. Suppose $AB\ge 0$ (positive semi-definite) and $Z>0$ with extremal eigenvalues $a$ and $b$. Then, the sharp inequality 
$$\Vert ZAB\Vert\le \,\frac{a+b}{2\sqrt{ab}}\,\Vert BZA\Vert
$$
holds for every unitarily invariant norm. Among the consequences, we get the operator inequality $XZX \le [(a+b)^2/4{ab}]Z$ for every $0\le X\le I$, and some Kantorovich type inequalities (Mond-Pe$\check{{\rm c}}$ari\'c inequalities). Also in connection, reverse inequalities of Davis and Hansen-Pedersen characterizations of operator convexity are established. For instance, given any  operator convex function $f:[0,\infty)\longrightarrow [0,\infty)$ and any subspace ${\cal E}$,
$$
f(Z_{\cal E}) \ge \,\frac{4ab}{(a+b)^2}\,(f(Z))_{\cal E}.
$$
In passing, we point out a simplified proof of Hansen-Pedersen's inequality.

Mathematical subjects classification:  15A60 47A30 47A63}

\vskip 25pt
{\large\bf \ Introduction}
\vskip 10pt
Capital letters $A$, $B\dots Z$ mean $n$-by-$n$ complex matrices, or operators on a finite dimensional Hilbert space ${\cal H}$; $I$ stands for the identity. When $A$ is positive semidefinite, resp.\ positive definite, we write $A\ge 0$, resp.\ $A>0$. Let $\Vert\cdot\Vert$ be a general symmetric (or unitarily invariant) norm, i.e.\ $\Vert UAV\Vert=\Vert A\Vert$ for all $A$ and all unitaries $U$, $V$. If $A$ and $B$ are such that the product $AB$ is normal, then a classical inequality claims [1, p.\ 253]
\begin{equation}
\Vert AB\Vert\le \Vert BA\Vert
\end{equation}

 Section 1  presents a generalization of (1) when $AB\ge0$. Then, for $Z>0$,
\begin{equation}
\Vert ZAB\Vert\le \,\frac{a+b}{2\sqrt{ab}}\,\Vert BZA\Vert
\end{equation}
where $a$, $b$ are the extremal eigenvalues of $Z$. Several sharp inequalities are derived.
For instance, if $0\le X\le I$, then
$$
XZX \le \,\frac{(a+b)^2}{4ab}\, Z.
$$
Another example concerns compressions $Z_{\cal E}$ of $Z$ onto subspaces ${\cal E}\subset{\cal H}$,
\begin{equation}
(Z_{\cal E})^{-1} \ge \,\frac{4ab}{(a+b)^2}\,(Z^{-1})_{\cal E}.
\end{equation}
\vskip 5pt
This  Kantorovich type inequality is due to Mond-Pe$\check{{\rm c}}$ari\'c. In Section 2 we extend (3) to all operator convex functions $f:[0,\infty)\longrightarrow [0,\infty)$. Such inequalities are reverse inequalities to Davis' characterization of operator convexity via compressions. Equivalently, we show that, given any isometric column of operators $\{A_i\}_{i=1}^m$, i.e.\ $\sum A_i^*A_i=I$, we have
$$
f(\sum A_i^*Z_iA_i) \ge \,\frac{4ab}{(a+b)^2}\,\sum A_i^*f(Z_i)A_i.
$$
This is a reverse inequality to the Hansen-Pedersen inequality.

\vskip 20pt\noindent
{\large\bf 1. Norms inequalities}

\vskip 10pt\noindent
{\bf Lemma 1.1.} {\it Let  $Z>0$ with extremal eigenvalues  $a$ and $b$. Then, for every norm one vector $h$,
$$
\Vert Zh\Vert \le \,\frac{a+b}{2\sqrt{ab}}\, \langle h,Zh\rangle.
$$
}
 
\vskip 10pt\noindent
{\bf Proof.} Let ${\cal E}$ be any subspace of ${\cal H}$ and let $a'$ and $b'$ be the extremal eigenvalues of $Z_{\cal E}$. Then $a\ge a'\ge b'\ge b$ and, setting $t=\sqrt{a/b}$, $t'=\sqrt{a'/b'}$, we have $t\ge t'\ge1$. Since $t\longrightarrow t+1/t$ increases on $[1,\infty)$ and
$$
\frac{a+b}{2\sqrt{ab}}\,=\frac{1}{2}\left( t+ \frac{1}{t}\right), \qquad
 \frac{a'+b'}{2\sqrt{a'b'}}\, =\frac{1}{2}\left( t'+ \frac{1}{t'}\right),
$$
we infer
$$
\frac{a+b}{2\sqrt{ab}}\,\ge \frac{a'+b'}{2\sqrt{a'b'}}\,.
$$
Therefore, it suffices to prove the lemma for $Z_{\cal E}$ with ${\cal E}={\rm span}\{h,Zh\}$. Hence, we may assume $\dim{\cal H}=2$, $Z=ae_1\otimes e_1+be_2\otimes e_2$ and $h=xe_1 +(\sqrt{1-x^2})e_2$. Setting
$x^2=y$ we have
$$
\frac{||Zh||}{\langle h, Zh\rangle}=\frac{\sqrt{a^2y+b^2(1-y)}}{ay+b(1-y)}.
$$
The right hand side attains its maximum on $[0,1]$ at $y=b/(a+b)$,
and then
$$
\frac{||Zh||}{\langle h, Zh\rangle}=\,\frac{a+b}{2\sqrt{ab}}
$$
proving the lemma.\qquad $\Box$

\vskip 10pt\noindent
{\bf Theorem 1.2.} {\it Let  $A$, $B$ such that $AB\ge0$. Let $Z>0$ with extremal eigenvalues  $a$ and $b$. Then, for every symmetric norm, the following sharp inequality holds
$$
\Vert ZAB\Vert \le \,\frac{a+b}{2\sqrt{ab}}\, \Vert BZA\Vert.
$$
}

\vskip 10pt\noindent
{\bf Proof.} For the sharpness see Remark 1.9 below.

It suffices to consider the Fan $k$-norms $\Vert\cdot\Vert_{(k)}$ [1, p.\ 93]. Fix $k$ and let $\Vert\cdot\Vert_1$ denote the trace-norm. There exist two rank $k$ projections $E$ and $F$ such that
\begin{align*} 
\Vert ZAB\Vert_{(k)} &=\Vert ZABE\Vert_1 \\
&=\Vert Z(AB)^{1/2}F(AB)^{1/2}E\Vert_1 \\
&\le \Vert Z(AB)^{1/2}F(AB)^{1/2}\Vert_1.
\end{align*}
Consider the canonical decomposition
$$
(AB)^{1/2}F(AB)^{1/2}=\sum_{j=1}^kc_j\, h_j\otimes h_j
$$
in which $\{h_j\}_{j=1}^k$ is an orthonormal system and $\{h_j\otimes h_j\}_{j=1}^k$ are the associated rank one projections. We have, using the triangle inequality and then the above lemma, 
\begin{align*}
\Vert Z(AB)^{1/2}F(AB)^{1/2}\Vert_1 &\le \sum_{j=1}^k c_j\Vert Z h_j\otimes h_j\Vert_1 \\
&= \sum_{j=1}^k c_j\Vert Z h_j\Vert  \\
 &\le \,\frac{a+b}{2\sqrt{ab}}\,\sum_{j=1}^k c_j\langle h_j,Zh_j\rangle \\
 &= \,\frac{a+b}{2\sqrt{ab}}\,{\rm Tr} \,(AB)^{1/2}F(AB)^{1/2}Z.
 \end{align*}
 Next, there exists a rank $k$ projection $G$
 such that 
\begin{align*}
\frac{a+b}{2\sqrt{ab}}\,{\rm Tr} \,(AB)^{1/2}F(AB)^{1/2}Z&=\,\frac{a+b}{2\sqrt{ab}}\,{\rm Tr} \,(AB)^{1/2}F(AB)^{1/2}ZG \\
&\le \,\frac{a+b}{2\sqrt{ab}}\,{\rm Tr} \,GZ^{1/2}ABZ^{1/2}G \\
&\le \,\frac{a+b}{2\sqrt{ab}}\,\Vert Z^{1/2}ABZ^{1/2}\Vert_{(k)} \\
&\le \,\frac{a+b}{2\sqrt{ab}}\,\Vert BZA\Vert_{(k)}
\end{align*}
where at the last step we used the basic inequality (1).\qquad $\Box$

\vskip 10pt
 One may ask wether our theorem can be improved to singular values inequalities. This is not possible as it is shown by the next example: 
 
\vskip 10pt\noindent
Take
$$
A=\begin{pmatrix} 1&0 \\ 0&4 \end{pmatrix}, \ B=\begin{pmatrix} 4&0 \\ 0&1 \end{pmatrix}, \ 
Z=\begin{pmatrix} 5&3 \\ 3&5\end{pmatrix}.
$$
Then the largest and smallest eigenvalues of $Z$ are $a=8$ and $b=2$, so 
$$
\frac{a+b}{2\sqrt{ab}}=1.25.
$$
Besides, $\mu_2(ZAB)=8$ and $\mu_2(AZB)=4.604$, and since $4.604\times 1.25=5.755<8$, Theorem 1.2 can not be extended to singular values inequalities.

\vskip 10pt
We denote by ${\rm Sing(X)}$ the sequence of the singular values of $X$, arranged in decreasing order and counted with their multiplicities. Similarily, when $X$ has only real eigenvalues, ${\rm Eig}(X)$ stands for the sequence of $X$'s eigenvalues. Given two sequences of real numbers $\{a_j\}_{j=1}^n$ and 
$\{b_j\}_{j=1}^n$, we use the notation $\{a_j\}_{j=1}^n\prec_w\{b_j\}_{j=1}^n$ for weak-majorisation, that is
$\sum_{j=1}^ka_j\le\sum_{j=1}^kb_j$, $k=1,\dots$.
\vskip 10pt
A straightforward application of Theorem 1.2 is:
\vskip 10pt\noindent
{\bf Corollary 1.3.} {\it Let $A\ge0$ and let $Z>0$ with extremal eigenvalues  $a$ and $b$. Then, 
$$
{\rm Sing}(AZ) \prec_w \,\frac{a+b}{2\sqrt{ab}} \,\,{\rm Eig}(AZ).
$$
}

\vskip 10pt\noindent
{\bf Proof.} For each Fan norms, replace $A$ and $B$ by $A^{1/2}$ in Theorem 1.2.\qquad $\Box$

\vskip 10pt 
Special cases of the above corollary are:

\vskip 10pt\noindent
{\bf Corollary 1.4.} {\it Let $A\ge0$ and let $Z>0$ with extremal eigenvalues  $a$ and $b$. Then, 
$$
\Vert AZ\Vert_{\infty}\le \,\frac{a+b}{2\sqrt{ab}}\,\rho(AZ) 
$$
and
$$
\Vert AZ\Vert_1\le \,\frac{a+b}{2\sqrt{ab}}\,{\rm Tr}\,AZ. 
$$
}

\vskip 10pt\noindent
Here, $\Vert\cdot\Vert_{\infty}$ stands for the standard operator norm and $\rho(\cdot)$ for the spectral radius.

\vskip 10pt
 From the preceding result, one may derive an interesting operator inequality:
\vskip 10pt\noindent
{\bf Corollary 1.5.} {\it Let $0\le A\le I$ and let $Z>0$ with extremal eigenvalues  $a$ and $b$. Then, 
$$
AZA \le \,\frac{(a+b)^2}{4ab}\, Z.
$$
}

\vskip 10pt\noindent
{\bf Proof.} The claim is equivalent to the operator norm inequalities
$$
\Vert Z^{-1/2}AZAZ^{-1/2}\Vert_{\infty} \le \,\frac{(a+b)^2}{4ab}
$$
or
$$
\Vert Z^{-1/2}AZ^{1/2}\Vert_{\infty} \le \,\frac{a+b}{2\sqrt{ab}}.
$$
But the previous corollary entails
\begin{align*}
\Vert Z^{-1/2}AZ^{1/2}\Vert_{\infty}&=\Vert Z^{-1/2}AZ^{-1/2}Z\Vert_{\infty} \\
&\le \,\frac{a+b}{2\sqrt{ab}}\,\rho(Z^{-1/2}AZ^{-1/2}Z) \\
& =\,\frac{a+b}{2\sqrt{ab}}\,\Vert A\Vert_{\infty} \\
& \le\,\frac{a+b}{2\sqrt{ab}},
\end{align*}
hence, the result holds. \qquad $\Box$

\vskip 10pt
A special case of Corollary 1.5 gives a comparison bewtween $Z$ and the compression $EZE$, for an arbitrary projection $E$.

\vskip 10pt\noindent
{\bf Corollary 1.6.} {\it Let $Z>0$ with extremal eigenvalues  $a$ and $b$ and let $E$ be any projection. Then, 
$$
EZE \le \,\frac{(a+b)^2}{4ab}\, Z.
$$
}

\vskip 10pt\noindent We may then derive a classical inequality:

\vskip 10pt\noindent
{\bf Corollary 1.7.} (Kantorovich) {\it Let $Z>0$ with extremal eigenvalues  $a$ and $b$ and let $h$ be any norm one vector. Then, 
$$
\langle h, Zh\rangle \langle h, Z^{-1}h\rangle \le \,\frac{(a+b)^2}{4ab}.
$$
}

\vskip 10pt\noindent
{\bf Proof.} Rephrase Corollary 1.6 as 
$$
\Vert Z^{-1/2}EZEZ^{-1/2}\Vert_{\infty}\le \,\frac{(a+b)^2}{4ab}
$$
and take $E=h\otimes h$. \qquad $\Box$

\vskip 10pt
A classical inequality in Matrix theory, for positive definite matrices, claims that "The inverse of a principal submatrix is less than or equal to the corresponding submatrix of the inverse" [6, p.\ 474]. In terms of compressions, this means
\begin{equation}
(Z_{\cal E})^{-1} \le (Z^{-1})_{\cal E}
\end{equation}
for every subspace ${\cal E}$ and every $Z>0$. Corollary 1.6 entails a reverse inequality, first proved by B.\ Mond and J.E.\ Pe$\check{{\rm c}}$ari\'c [7]:

\vskip 10pt\noindent
{\bf Corollary 1.8.}  (Mond-Pe$\check{{\rm c}}$ari\'c) {\it Let $Z>0$ with extremal eigenvalues  $a$ and $b$. Then, for every subspace ${\cal E}$, 
$$
(Z_{\cal E})^{-1} \ge \,\frac{4ab}{(a+b)^2}\,(Z^{-1})_{\cal E}.
$$
}

\vskip 10pt\noindent
Note that Corollary 1.8 implies Corollary 1.7.
\vskip 10pt\noindent
{\bf Proof.} Let $E$ be the projection onto ${\cal E}$. By Corollary 1.6, for every $r>0$, there exists $x>0$ such that
$$
EZE + x E^{\perp} \le \,\frac{(a+b)^2}{4ab}\, (Z + rI).
$$
Since $t\longrightarrow -1/t$ is operator monotone we deduce
$$
(EZE + x E^{\perp})^{-1} \ge \,\frac{4ab}{(a+b)^2}\, (Z + rI)^{-1}
$$
so that
$$
(Z_{\cal E})^{-1} \ge \,\frac{4ab}{(a+b)^2}\, \{(Z + rI)^{-1}\}_{\cal E}
$$
and the result follows by letting $r\longrightarrow 0$. \qquad $\Box$

\vskip 10pt\noindent
{\bf Remark 1.9.} All the previous inequalities  are sharp. Indeed, let $h$ be a norm one vector for which equality occurs in Lemma 1.1. Then, replacing A, B, E by $h\otimes h$ and ${\cal E}$ by 
${\rm span}\{h\}$ in the above statements, yields equality cases. 

\vskip 10pt\noindent
{\bf Remark 1.10.} As for a standard proof of (1) [1, p.\ 253], it is tempting to first prove Theorem 1.2 for the operator norm and then to use an antisymmetric tensor product argument to derive the general case. Such an approach seems impossible. Indeed if $a_k$ and $b_k$ are the extremal eigenvalues of $\wedge^k(Z)$, then the relation
$$
\frac{(a_k+b_k)^2}{4a_kb_k}\,\le \,\left(\frac{(a+b)^2}{4ab}\right)^{k}
$$
is not true in general.

\vskip 10pt
The next result states a companion inequality  to Corollary 1.8.

\vskip 10pt\noindent
{\bf Proposition 1.11.} {\it Let $Z>0$ with extremal eigenvalues  $a$ and $b$ and let $1\le p\le 2$. Then, for every subspace ${\cal E}$, 
$$
(Z_{\cal E})^p \ge \,\frac{4ab}{(a+b)^2}\,(Z^p)_{\cal E}.
$$
}

\vskip 10pt\noindent
{\bf Proof.} Let $E$ be the projection onto ${\cal E}$. For any norm one vector $h\in{\cal E}$, Lemma 1.1
implies
\begin{align*}
\langle h, (Z^p)_{\cal E}h\rangle &= \langle h, EZ^pEh\rangle \\
&= \Vert Z^{p/2}h\Vert^2 \\ 
&\le \,\frac{(a+b)^2}{4ab}\, \langle h, Z^{p/2}h\rangle^2. 
\end{align*}
Then, using the concavity of $t\longrightarrow t^{p/2}$ and next the convexity of $t\longrightarrow t^p$, we deduce 
\begin{align*}
\langle h, (Z^p)_{\cal E}h\rangle &\le \,\frac{(a+b)^2}{4ab}\,  \langle h, Zh\rangle^p \\
&=\frac{(a+b)^2}{4ab}\,  \langle h,EZEh\rangle^p \\ 
&\le \,\frac{(a+b)^2}{4ab}\, \langle h, (Z_{\cal E})^ph\rangle. 
\end{align*}
and the proof is complete.\qquad $\Box$

\vskip 20pt\noindent
{\large\bf 2. Operator convexity}

\vskip 10pt
 Davis' characterization of operator convexity [2] claims: $f$ is operator convex on $[a,b]$ if and only if for every subspace ${\cal E}$ and every Hermitian $Z$ with spectrum in $[a,b]$, 
\begin{equation*}
f(Z_{\cal E}) \le (f(Z))_{\cal E} \tag{D}
\end{equation*}
Since $t\longrightarrow t^p$, $1\le p\le 2$ and $t\longrightarrow 1/t$ are operator convex on $(0,\infty)$, both Proposition 1.11 and Corollary 1.8 are reverse inequalities to Davis' characterization of operator convexity.

\vskip 10pt
Proposition 1.11 is a special case of the next theorem.

\vskip 10pt\noindent
{\bf Theorem 2.1.} {\it Let $f:[0,\infty)\longrightarrow[0,\infty)$ be operator convex and let $Z>0$ with extremal eigenvalues $a$ and $b$. Then, for every subspace ${\cal E}$, 
$$
f(Z_{\cal E}) \ge \,\frac{4ab}{(a+b)^2}\,(f(Z))_{\cal E}.
$$
}

\vskip 10pt\noindent
{\bf Proof.} We have the integral representation [1]
$$
f(t)=\alpha+\beta t+\gamma t^2 +\int_0^{\infty}\frac{\lambda t^2}{\lambda +t}\,d\mu(\lambda),
$$
where $\alpha,\beta,\gamma$ are nonnegative scalars and $\mu$ is a positive finite measure. Therefore,
it suffices to prove the result for 
$$
\alpha+\beta t+\gamma t^2 
$$
and
$$
f_{\lambda}(t)=\frac{\lambda t^2}{\lambda +t}.
$$
The quadratic case is a staightforward application of Proposition 1.11. To prove the $f_{\lambda}$ case, note that $f_{\lambda}$ is convex meanwhile $f^{1/2}_{\lambda}$ is convave and then proceed as in the proof of Proposition 1.11. \qquad $\Box$
 
\vskip 10pt
Davis' characterization (D) of operator convexity  is equivalent to the following result of Hansen-Pedersen [5]. 

Recall that a family $\{ A_i\}_{i=1}^m$ form an isometric column when $\sum A^*_iA_i=I$. 

\vskip 10pt\noindent
{\bf Theorem 2.2.} (Hansen-Pedersen) \ {\it Let $\{ Z_i\}_{i=1}^m$ be Hermitians with spectrum lying in $[a,b]$ and let $f$ be operator convex  $[a,b]$. Then, for every isometric column $\{ A_i\}_{i=1}^m$, 
\begin{equation*}
f(\sum A_i^*Z_iA_i) \le \sum A_i^*f(Z_i)A_i. \tag{Jo}
\end{equation*}
}
\vskip 10pt
 (Jo)  is the operator version of Jensen's inequality: operator convex combinations and operator convex functions replace the ordinary ones. As a sthraightforward consequence, we have the following contractive version of (Jo):
 
\newpage
    
\vskip 10pt\noindent
{\bf Corollary 2.3.} (Hansen-Pedersen) \ {\it Let $\{ Z_i\}_{i=1}^m$ be Hermitians with spectrum lying in $[a,b]$ and let $f$ be operator convex  $[a,b]$ with $0\in[a,b]$ and $f(0)\le0$. Then, for every contraction $A$, 
\begin{equation*}
f(A^*ZA) \le  A^*f(Z)A. \tag{C}
\end{equation*}
}

Exactly as Theorem 2.1 is a reverse inequality to (D),  the following  results is a reverse inequality to (Jo).

\vskip 10pt\noindent
{\bf Theorem 2.4.} {\it Let $f:[0,\infty)\longrightarrow[0,\infty)$ be operator convex and let $\{Z_i\}_{i=1}^m$ be positive with spectrum lying in  $[a,b]$, $a>0$. Then, for every isometric column $\{A_i\}_{i=1}^m$, 
$$
f(\sum A_i^*Z_iA_i) \ge \,\frac{4ab}{(a+b)^2}\,\sum A_i^*f(Z_i)A_i.
$$
}

\vskip 10pt\noindent

Let us consider a very special case: For every $A,B>0$ with spectrum lying on $[r,2r]$, $r>0$, and for every
operator convex  $f:[0,\infty)\longrightarrow[0,\infty)$, we have
$$
\frac{8}{9} \cdot \frac{f(A)+f(B)}{2}\le f\left(\frac{A+B}{2}\right)\le \frac{f(A)+f(B)}{2}.
$$
The left inequality gives a negative answer to an approximation problem: Let $f$ be an operator convex function on $[a,b]$, $0<a<b$, and let $\varepsilon>0$. Then, in general, there is no operator convex function $g$ on [$0,\infty)$ such that
$$
\max_{x\in[a,b]} |f(x)-g(x)|< \varepsilon.
$$
\vskip 10pt
From Theorem 2.4 we obtain a reverse inequality to (C):
\vskip 10pt\noindent
{\bf Corollary 2.5.} {\it Let $f:[0,\infty)\longrightarrow[0,\infty)$ be operator convex and let $Z>0$ with extremal eigenvalues $a$ and $b$. Then, for every contraction $A$, 
$$
f(A^*ZA) \ge \,\frac{4ab}{(a+b)^2}\, A^*f(Z)A.
$$
}

\vskip 10pt
We turn to the proof of Theorem 2.4 and Corollary 2.5.
\vskip 10pt\noindent
{\bf Proof.} Consider the following operators acting  on $\oplus^m{\cal H}$,
$$V=
\begin{pmatrix}
A_1&0&\cdots&0 \\ \vdots &\vdots &\ &\vdots \\
A_m &0 &\cdots &0
\end{pmatrix},\qquad \tilde{Z}=
\begin{pmatrix}
Z_1&\ &\  \\ 
\ &\ddots &\  \\
\ &\  &Z_m
\end{pmatrix} 
$$
and note that $V$ is a partial isometry.
Denoting by ${\cal H}$ the first summand of the direct sum $\oplus^m{\cal H}$ and by $X\!:\!{\cal H}$ the restriction of $X$ to ${\cal H}$, we observe that
$$
f(\sum A_i^*Z_iA_i)= f(V^*\tilde{Z}V)\!:\!{\cal H} = V^*f(\tilde{Z}_{V({\cal H})})V\!:\!{\cal H}.
$$
Applying  Theorem 2.1 with ${\cal E}=V({\cal H})$, we get
\begin{align*}
f(\sum A_i^*Z_iA_i) &\ge \,\frac{4ab}{(a+b)^2}\, V^*f(\tilde{Z})_{V({\cal H})}V\!:\!{\cal H} \\
&=\,\frac{4ab}{(a+b)^2}\,\sum A_i^*f(Z_i)A_i.
\end{align*}
and the proof of Theorem 2.4 is complete. 
To obtain its corollary, take an operator $B$ such that $A^*A+B^*B=I$. Then, note that, using $f(0)\ge0$, 
\begin{align*}
f(A^*ZA)=f(A^*ZA+B^*0B) &\ge \,\frac{4ab}{(a+b)^2}\, \{A^*f(Z)A +B^*f(0)B\} \\
&\ge \,\frac{4ab}{(a+b)^2}\, A^*f(Z)A 
\end{align*}
by application of Theorem 2.4. \qquad $\Box$

\vskip 10pt\noindent
{\bf Remark 2.6.} Corollary 1.8 and Proposition 1.11 for $p=2$ have been obtained by Mond-Pecaric in the more general form of Theorem 2.4. Note that Proposition 1.11 with $p=2$ immediately implies Lemma 1.1; hence, we have no pretention of originality in establishing this basic lemma.

\vskip 10pt\noindent
{\bf Remark 2.7.} Hansen-Pedersen first prove the contractive version (C) in [4] and then, some twenty years later [5], prove the more general form (Jo). When proving (C) they noted a technical difficulty to derive (Jo) when $0\notin [a,b]]$. In fact, this difficulty can be easily overcomed: Note that if (Jo) is valid for every operator convex functions on an interval $[a,b]$, then (Jo) is also valid on every interval of the type $[a+r,b+r]$.

\vskip 10pt\noindent
{\bf Remark 2.8.} (D), (Jo), (C) are equivalent statements. Similarly, Theorems 2.1, 2.4 	and Corollary 2.5 are equivalent.

\vskip 10pt
Clearly, the previous results can be suitably restated for operators acting on infinite dimensional spaces.

\vskip 10pt
Inspired by the seminal paper [3], we note that Corollary 2.5 can be stated in a still more general framework. Let ${\rm B}({\cal H})$ denote the algebra of all (bounded) linear operators on a  separable Hilbert space ${\cal H}$.

\newpage
\vskip 10pt\noindent
{\bf Corollary 2.9.} {\it Let $\Phi:{\cal Z}\longrightarrow{\rm B}({\cal H})$ be a positive, linear contraction on a $C^*$-algebra ${\cal Z}$. Let $Z\in{\cal Z}$, $Z>0$ with ${\rm Sp}(Z)\subset[a,b]$, $a>0$. Then, for every operator convex function $f:[0,\infty)\longrightarrow[0,\infty)$, 
$$
f\circ\Phi(Z) \ge \,\frac{4ab}{(a+b)^2}\, \Phi\circ f(Z).
$$
}

\vskip 10pt\noindent
{\bf Proof.} Restricting $\Phi$ to the commutative $C^*$-subalgebra generated by $Z$, one may suppose $\Phi$ completely positive. By Stinepring's dilation Theorem [8], there exist a larger Hilbert space ${\cal F}\supset{\cal H}$, a linear contraction $A:{\cal H}\longrightarrow{\cal F}$ and a $*$-homomorphism $\pi:{\cal Z}\longrightarrow{\rm B}({\cal F})$ such that $\Phi(\cdot)=A^*(\pi(\cdot))_{\cal F}A$. Therefore
\begin{align*}
f\circ\Phi(Z) &= f(A^*\pi(Z)A) \\
&\ge \,\frac{4ab}{(a+b)^2}\, A^*f(\pi(Z))A \\
&=\,\frac{4ab}{(a+b)^2}\, A^*\pi(f(Z))A \\
&=\,\frac{4ab}{(a+b)^2}\, \Phi\circ f(Z)
\end{align*}
where at the second step we apply Corollary 2.5 which can be extended to this situation by inspection of the proof of Theorem 2.4 and Corollary 2.5. \qquad $\Box$

\vskip 15pt
{\bf References}
\vskip 5pt
\noindent
{\small 
\noindent
[1] R.\ Bhatia, Matrix Analysis, Springer, Germany, 1996.

\noindent
[2] C.\ Davis, A Shwarz inequality for convex operator functions, Proc.\ Amer.\ Math.\ Soc.\ 8 (1957) 42-44.

\noindent
[3] F.\ Hansen, An operator inequality, Math.\ Ann.\ 246 (1980) 249-259.

\noindent
[4] F.\ Hansen and G.\ K.\ Pedersen, Jensen's inequality for operator sand Lowner's Theorem, Math. Ann. 258 (1982) 229-241.

\noindent
[5] F.\ Hansen and G.\ K.\ Pedersen, Jensen's operator inequality, Bull.\ London Math.\ Soc.\ 35 (2003) 553-564.

\noindent
[6] R.A.\ Horn,  C.R.\ Johnson, Matrix Analysis, Cambridge Univ.\ Press, Cambridge, 1985.

\noindent
[7] B.\ Mond,  J.E.\ Pe$\check{{\rm c}}$ari\'c, A matrix version of the Ky Fan generalization of the Kantorovich inequality,
Linear and Multilinear Algebra 36 (1994) 217-221.

\noindent
[8] W.\ F.\ Stinepring, Positive functions on $C^*$-algebras, Proc.\ Amer.\ Math.\ Soc.\ 6 (1955) 211-216.  }

\vskip 15pt
\begin{flushright}
{\it 
 Jean-Christophe Bourin
\vskip 10pt
 Universit\'e de Cergy-Pontoise, d\'ept.\ de Math\'ematiques
\vskip 5pt
 2 rue Adolphe Chauvin, 95302 Pontoise
}
\end{flushright}

\end{document}